\documentclass[twoside,
]{amsart}

  \usepackage[utf8]{inputenc}
  \usepackage{amsmath}
  \usepackage{amssymb}
  \usepackage{amsthm}
  \usepackage{mathrsfs}
  \usepackage{bbm}
  \usepackage{pifont}
  \usepackage{graphicx}
  \usepackage{xmpmulti}
  \usepackage{color}

	\newtheoremstyle{slanted}
	{}
	{}
	{\slshape}
	{}
	{\bfseries}
	{.}
	{ }
	{}
	
	\theoremstyle{slanted}
	\newtheorem{theo}{Theorem}[section]

	\newtheorem{lemma}[theo]{Lemma}

	\def\ind#1{\mathbbmss{1}_{#1}}

	\newcommand{\EE}{\mathbb{E}}

	\def\ind#1{\mathbbmss{1}_{#1}}
	\newcommand{\ZZ}{\mathbb{Z}}
	\renewcommand{\AA}{\mathbb{A}}
	\newcommand{\BB}{\mathbb{B}}
	
	\newcommand{\RR}{\mathbb{R}}
	\newcommand{\PP}{\mathbb{P}}
	\newcommand{\NN}{\mathbb{N}}
	
	\newcommand{\var}{\mathop{\mbox{var}}}
	\newcommand{\dbar}{\bar d}
	\newcommand{\xiell}{\left(\xi\wedge\ell\right)}
	\newcommand{\zetaell}{\left(\zeta\wedge\ell\right)}
\newcommand{\hkr}{h_{\mbox{\scriptsize Kr}}} 

\newcommand{\n}[1]{{N}_{\!#1}}

\title{Zero Krengel Entropy does not kill Poisson Entropy}

\author{Élise Janvresse  and Thierry de la Rue}

\address{Laboratoire de Math\'ematiques Rapha\"el Salem\\
	UMR 6085 CNRS -- Universit\'e de Rouen\\
	Avenue de l'Universit\'e\\
	B.P. 12\\
	F76801 Saint-\'Etienne-du-Rouvray Cedex}
\email{Elise.Janvresse@univ-rouen.fr}  \email{Thierry.de-la-Rue@univ-rouen.fr}

\begin{document}
\bibliographystyle{amsplain}

\begin{abstract}
We prove that the notions of Krengel entropy and Poisson entropy for infinite-measure-preserving transformations do not always coincide: We construct a conservative infinite-measure-preserving transformation with zero Krengel entropy (the induced transformation on a set of measure 1 is the Von Neumann-Kakutani odometer), but whose associated Poisson suspension has positive entropy.
\end{abstract}

\maketitle

\section{Introduction}
\subsection{Entropy for infinite-measure-preserving transformations}

There exist several notions of entropy for infinite transformations, which elegantly generalize Kolmogorov's entropy of a probability-preserving transformation.
Krengel~\cite{krengel67} comes down to the finite-measure case by considering the entropy of the induced transformation on a set of finite measure:
The \emph{Krengel entropy} of a conservative measure-preserving
transformation $(X,\mathcal{B},\mu,T)$ is defined as
$$
\hkr (X,\mathcal{B},\mu,T) :=\sup_{A \in\mathcal{F}_+}\mu(A)\,h(A,\mathcal{B}\cap A,{\mu}_A,T_A),
$$
where $\mathcal{F}_+$ is the collection of sets in $\mathcal{B}$ with
finite positive measure, ${\mu}_A$ is the normalized
probability measure on $A$ obtained by restricting $\mu$ to
$\mathcal{B} \cap A$, and $T_A:A \to A$ is the induced map on
$A$. Recall that this map is defined by
$$T_A(x) := T^{r_A(x)}(x),$$
where $r_A(x) := \min\{k \ge 1:~ T^k(x) \in A\}$ is the \emph{first-return-time map} associated to $A$.
As soon as $T$ is not purely periodic, Krengel proved that
$$
\hkr (X,\mathcal{B},\mu,T) =\mu(A)\,h(A,\mathcal{B}\cap A,\mu_A,T_A),
$$
where $A$ is any finite-measure \emph{sweep-out} set (\textit{i.e.} a set such that
${\bigcup_{n\ge0}}T^{-n}A=X$).

The \emph{Parry entropy} of an infinite-measure-preserving transformation $T$ has been defined in~\cite{parry_generators69} as the supremum of the conditional entropy of $\mathcal C$ with respect to $T^{-1}\mathcal C$, for all $\sigma$-finite sub-$\sigma$-algebras $\mathcal C$ such that $T^{-1}\mathcal C\subset \mathcal C$.

Recall now that to each infinite-measure-preserving transformation $T$ we can associate a probability-preserving transformation $T_*$ called its \emph{Poisson suspension}, and which can be described as follows (we refer to~\cite{roy_thesis} for details): We consider a Poisson process on $X$ with intensity $\mu$, which we can consider as a random collection of particles. These particles are distributed over $X$ in such a way that, denoting by $\n{B}$ the random number of particles in any finite-measure set $B$, for any finite collection of pairwise disjoint, finite-measure sets $B_1,\ldots,B_n\in\mathcal{B}$, the random variables $\n{B_1},\ldots,\n{B_n}$ are independent, and follow Poisson distributions with respective parameters $\mu(B_1),\ldots,\mu(B_n)$.  Then $T_*$ is defined on the canonical space of this Poisson process, and it consists in moving individually each of these particle according to the transformation $T$ on $X$.  The \emph{Poisson entropy} of an infinite-measure-preserving transformation was defined by Roy~\cite{roy_thesis} as the Kolmogorov entropy of its Poisson suspension.

Relations between these notions of entropy are studied in~\cite{jmrr08}:
On large classes of transformations (\textit{e.g.} quasi-finite transformations, rank-one transformations), it is proved that Poisson entropy is equal to Krengel entropy and to Parry entropy. Moreover, in any case, Parry entropy is dominated by both Krengel and Poisson entropy. 

It was asked in~\cite{jmrr08} whether, for any conservative measure-preserving transformation, these three definitions always coincide. The purpose of the present paper is to show that the answer is negative, by constructing a counterexample. 

\begin{theo}
 There exists a conservative infinite-measure-preserving transformation with zero Krengel entropy (hence zero Parry entropy), but whose associated Poisson suspension has positive entropy.
\end{theo}

%
%

\section{Construction}
\label{Sec:construction}
\subsection{Von Neumann-Kakutani odometer}
\begin{figure}
 \input{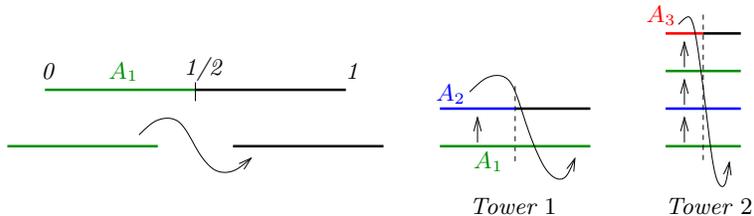}
\caption{First steps in the construction of the Von Neumann-Kakutani odometer by cutting and stacking.}
 \label{Fig:odometer}
\end{figure}

The transformation $T$ is constructed as a tower over the Von Neumann-Kakutani odometer $S$. 
Let us recall the construction of the latter by cutting and stacking (see Figure~\ref{Fig:odometer}). 
We start with the interval $A:=[0,1]$. 
The first step consists in cutting $A$ into two sub-intervals $A_1$ and $A\setminus A_1$ of measure $1/2$, and stacking $A\setminus A_1$ over $A_1$. 
We get a tower of height 2 which we call Tower 1. 
After step $n$, $A$ has been cut into $2^n$ sub-intervals which are stacked to get Tower $n$. 
This means that at this step each point of $A$ (except those lying on the top of the tower) is mapped by $S$ to the point of $A$ lying above it. 
We construct Tower $n+1$ by cutting Tower $n$ into two equal parts. 
We call $A_{n+1}$ the left half of the top interval of Tower $n$ and we stack the right part of the tower over $A_{n+1}$, thus dividing by 2 the measure of the set where $S$ is not yet defined.
Repeating this procedure defines the Von Neumann-Kakutani odometer which preserves the Lebesgue measure on $A$. 
It is well known that the odometer $S$ is ergodic and has zero entropy.

\subsection{Construction of $T$}
$T$ is constructed on $\RR_+$ in such a way that the induced transformation $T_A$ coincides with the odometer previously defined.
$T$ is completely defined (up to isomorphism) by giving for each point $x\in A$ the first return time $r_A(x)$ to $A$.

We fix an increasing sequence of integers $(M_n)_{n\ge0}$, with $M_n\to\infty$.
For any $n\ge1$, we choose a large enough integer $k_n$ (to be precised later). 
We define the first return time to $A$ so that its restriction to $A_n$ is uniformly distributed on $\{M_n, M_n+1, \ldots, M_n+k_n-1\}$ for any $n$.

We will see in section~\ref{Sec:positiveEntropy} that by choosing $k_n$ large enough, the entropy of $T_*$ is positive. 

\section{Poisson approximation lemma}

The purpose of this section is to prove the key lemma. This lemma roughly states that when $k_n$ is large enough, it is almost impossible in the Poisson suspension to keep track individually of the particles when they leave $A_n$ if we only have access to the number of particles in $A$. For this, we compare two processes: The first one is simply an i.i.d. sequence of Poisson random variables, whereas the second one modelizes particles leaving $A_n$ and coming back to $A$. The comparison between the two processes uses the notion of \emph{$\dbar$-distance}, of which we recall some properties.

\subsection{The $\dbar$-distance}
The $\dbar$-distance between two stationary processes has been introduced by Ornstein for the proof of the isomorphism theorem of Bernoulli shifts. We refer to~\cite{ornstein1974} or \cite{shields1996} for the  properties of this distance which we use later and which we recall here. 

Let $\xi$ and $\zeta$ be two stationary processes taking values in a countable alphabet $\BB$. For any integers $p<q$, we denote by $\xi|_p^q$ the finite sequence $(\xi_p,\xi_{p+1},\ldots,\xi_{q})$.

For $L\ge1$, let $J_L(\xi,\zeta)$ be the set of all joinings of $\xi|_1^L$ and $\zeta|_1^L$, that is probability distributions on $\BB^L\times\BB^L$ whose marginals are the distributions of $\xi|_1^L$ and $\zeta|_1^L$.
We first define $\dbar_L$ for any $L\ge1$, by
$$ 
\dbar_L (\xi,\zeta) := \min_{\lambda \in J_L(\xi,\zeta)}\EE_\lambda\Bigl[d_L\left(\xi|_1^L,\zeta|_1^L\right)\Bigr], 
$$
where $d_L$ is the Hamming distance between sequences of length $L$: 
$$ d_L(x_1\ldots x_L,z_1\ldots z_L):=\dfrac{1}{L}\sum_1^L \ind{x_i\neq z_i}. $$
Then, the $\dbar$-distance between $\xi$ and $\zeta$ is defined by
$$ \dbar (\xi, \zeta) := \sup_{L\ge1} \dbar_L (\xi,\zeta). $$
It can be shown that $\dbar (\xi, \zeta)$ is also the minimum of $\lambda\left(\xi_0\neq \zeta_0\right)$ when $\lambda$ ranges over all stationary joinings of $\xi$ and $\zeta$.

The two key properties of the $\dbar$-distance that we shall use are:
On the one hand, the fact that entropy of processes close in  $\dbar$-distance can be compared (Lemma~\ref{lemma:entropy} below). On the other hand, a practical tool to estimate the $\dbar$-distance between processes using conditional distributions on the past: If, for all large enough $n$, 
\begin{equation}
 \label{eq:conditional-property}
\sum_{b\in\BB} \left|\PP\Bigl( \xi_0=b \Big| \xi|_{-n}^{-1}\Bigr) - \PP\Bigl( \zeta_0=b \Big| \zeta|_{-n}^{-1}\Bigr)\right| < \varepsilon
\end{equation}
for all past $\xi|_{-n}^{-1}$ outside a set of measure $\varepsilon$ and all past $\zeta|_{-n}^{-1}$ outside a set of measure $\varepsilon$, then $\dbar(\xi,\zeta)\le 3 \varepsilon$. Moreover, the same conclusion holds if we replace in~\eqref{eq:conditional-property} the conditional distributions with respect to the past by conditional distributions with respect to finer $\sigma$-algebras.

\begin{lemma}
 \label{lemma:entropy}
Let $\xi$ be a stationary process taking values in a countable alphabet $\BB$ which has finite entropy. 
For any $\varepsilon>0$, there exists $\delta>0$ such that any stationary process $\zeta$ taking values in the same alphabet $\BB$ with $\dbar(\xi,\zeta)<\delta$ satisfies $h(\zeta)>h(\xi)-\varepsilon$.
\end{lemma}

\begin{proof}
 Without loss of generality we can assume $\BB=\ZZ^+$. For any integer $\ell>0$ we define the process $\xiell$ taking values in the finite alphabet $\{0,1,\ldots,\ell\}$ by 
$$ \xiell_x
=\begin{cases}
    \xi_x & \mbox{ if }\xi_x<\ell\\
    \ell & \mbox{ otherwise. }
   \end{cases}
$$
We choose $\ell$ large enough so that $h\xiell>h(\xi)-\varepsilon/2$.
Then we use the fact that entropy is a continuous function of processes taking values in a given finite alphabet, when these processes are topologized with the $\dbar$-distance (see \textit{e.g.} \cite{shields1996} page 100). Therefore we can find $\delta>0$ such that any process taking values in $\{0,1,\ldots,\ell\}$ at $\dbar$-distance at most $\delta$ from $\xiell$ has entropy at least $h(\xi)-\varepsilon$. Now, if $\dbar(\xi,\zeta)<\delta$, then $\dbar\left(\xiell,\zetaell\right)<\delta$ (where $\zetaell$ is defined from $\zeta$ in a similar way), hence 
$h(\zeta)\ge h\zetaell \ge h(\xi)-\varepsilon$.
\end{proof}

\subsection{Comparison between connected and disconnected processes}
Let $\AA$ be a finite alphabet, $P^B$ and $P^W$ be two probability measures on $\AA$ and $\lambda_{\AA\times\AA}$ be a joining of $P^B$ and $P^W$. Let $\delta$ be some fixed positive real number.
We define two processes $\xi$ and $\zeta$ on $\left(\NN^\AA\right)^\ZZ$.

The process $\xi$ is constructed from two independent sequences of i.i.d. random variables distributed according to the Poisson distribution of parameter $\delta/2$, which can be interpreted as numbers of black and white particles lying on each site of $\ZZ$. 
Then to each black (respectively white) particle we randomly and independently associate a label picked in $\AA$ according to $P^B$ (respectively $P^W$). 
For any $x\in\ZZ$ and any labels $a,b\in\AA$, $\xi^B_x(a)$  (respectively $\xi^W_x(b)$) is the total number of black (respectively white) particles labelled by $a$  (respectively $b$) at position $x$. 
In other words, the process $\xi$ associates in an i.i.d. way to each site $x\in\ZZ$ a finite sequence $\xi_x=\bigl(\xi^B_x(a),\xi^W_x(b)\bigr))_{a,b\in\AA}$ of independent random variables respectively distributed according to the Poisson distribution of parameter $P^B(a)\delta/2$ and $P^W(b)\delta/2$.

Let $k$ and $M$ be two integers. The process $\zeta$ is also constructed from black and white particles on $\ZZ$, but which are no longer independent. 
The number of black particles at each site is given by a sequence of i.i.d. random variables distributed according to the Poisson distribution of parameter $\delta/2$. 
For each black particle at position $x\in\ZZ$, we first pick a random integer $j$, uniformly in $\{M, M+1, \ldots, M+k-1\}$ and independently of all other particles. Then we link the black particle to a white particle that we put at position $x+j$. For each such couple of black and white particles, a couple of labels is picked in $\AA\times\AA$ with probability $\lambda_{\AA\times\AA}$: The first label is associated to the black particle and the second label to the white one. Then, for any $x\in\ZZ$, $\zeta_x=\bigl(\zeta^B_x(a),\zeta^W_x(b)\bigr))_{a,b\in\AA}$ denote the number of  black particle labelled by $a$ and the number of white particles labelled by $b$ at position $x$. 

\begin{lemma}
 \label{lemma:general}
For any $\varepsilon >0$, if $k$ is large enough, $\dbar (\xi, \zeta)<\varepsilon$.
\end{lemma}

\begin{proof}[Proof of Lemma~\ref{lemma:general}, simple case]
We first prove the lemma in the case where $\AA$ is reduced to a singleton. Since all particles have the same label, we just forget it and simply count the number of black and white particles on each site. 

We now prove that if $k$ is large enough, for any $n\ge M+k$,
$$
\sum_{j,\ell\in\NN} \left| P\left(\zeta^B_0=j,\ \zeta^W_0=\ell \ \big| \ \zeta |_{-n}^{-1} \right) - \frac{e^{-\delta/2}(\delta/2)^{j}}{j!}\frac{e^{-\delta/2}(\delta/2)^{\ell}}{\ell!}\right| <\varepsilon ,
$$
with probability $1-\varepsilon$ on $\zeta |_{-n}^{-1}$. 
Since $\zeta^B_0$ is Poisson distributed with parameter $\delta/2$, and independent from $\zeta |_{-n}^{-1}$, and since $\zeta^W_0$, is independent from $\zeta^B_0$ conditionally to the past, it is enough to prove that if $k$ is large enough, for any $n\ge M+k$,
$$
\sum_{\ell\in\NN} \left| P\left(\zeta^W_0=\ell \ \big| \ \zeta |_{-n}^{-1} \right) - \frac{e^{-\delta/2}(\delta/2)^{\ell}}{\ell!}\right| <\varepsilon ,
$$
with probability $1-\varepsilon$ on $\zeta |_{-n}^{-1}$. 

In fact we rather condition with respect to an enriched past (see Figure~\ref{Fig:enriched_past}): Assume that besides the number of black and white particles on each site $x\in \{-n, \ldots,-1\}$, we also know 
which the links are between them.

\begin{figure}[h]
\label{Fig:enriched_past}
\input{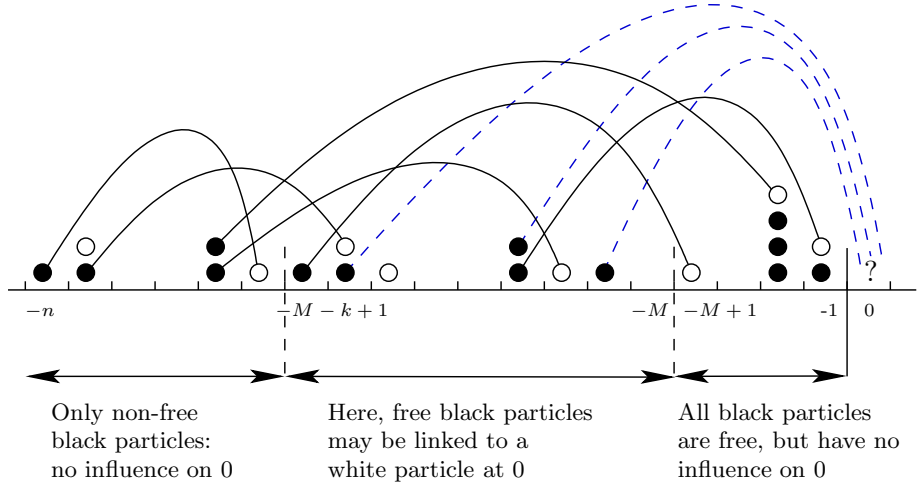}
\caption{The enriched past.}
\end{figure}

{From} what we know, we can distinguish two kinds of black particles between $-n$ and $-1$: Those which are linked to a white particle lying on the left of 0, and those, called \textit{free particles}, whose white particle's position is unknown.
Observe that only free particles may have some influence on $\zeta^W_0$. 
Hence, black particles lying on the left of site $-(M+k-1)$ have no influence on $\zeta^W_0$ since they are not free. 
Black particles lying on the right of site $-M$ are free but nevertheless have no influence on $\zeta^W_0$ since a black particle is linked to a white particle at distance at least $M$. 

So it remains to study the influence on $\zeta^W_0$ of free black particles at sites between $-(M+k-1)$ and $-M$. 
For $1\le j\le k$ let us denote by $F_j$ the number of free particles at site $-(M+k)+j$.
Any black particle at site $-(M+k)+j$ has probability $j/k$ to be free. 
Therefore, $F_j$ follows the Poisson distribution with parameter $\frac{\delta j}{2k}$.

Fix $1\le j\le k$. Assume there is a free particle at site $-(M+k)+j$. Since there are $j$ possible positions for its white particle, the latter has probability $1/j$ to lie at site 0.
Hence, conditionally to our enriched past, the number of white particles at site 0 can be written as 
$$
\sum_{j=1}^k \sum_{\ell=1}^{F_j} B_\ell^j,
$$
where $(B_\ell^j)$ are independent Bernoulli random variables with respective parameter $1/j$.
The law of such a sum of independent Bernoulli variables is close to a Poisson distribution of parameter $\delta/2$ as soon as the sum of the parameters is close to $\delta/2$ and all parameters are small enough (see~\cite{Billingsley86}, Theorem~23.2 page 312). Therefore we can choose a large enough integer $J$, and $\varepsilon_1$ small enough so that, if $Z$ is a sum of independent Bernoulli random variables, each with parameter less than $1/J$ and such that the sum of parameters is within $\varepsilon_1$ of $\delta/2$, then
$$ \sum_{\ell\ge0} \left| \PP(Z=\ell)-\exp(-\delta/2)\frac{(\delta/2)^\ell}{\ell !}\right| < \varepsilon. $$

We have now to avoid bad configurations, that is configurations of the enriched past which have free particles close to $-(M+k-1)$ giving rise to Bernoulli with large parameters, and configurations such that the sum of the parameters is not close enough to $\delta/2$. 

\subsubsection*{Control of the parameters' size}
We compute the probability that no free particle lie between $-(M+k)+1$ and $-(M+k)+J$: 
$$
\PP(F_j=0,\ j=1,\ldots,J) = \prod_{j=1}^{J} \exp\left(-\frac{\delta j}{2k}\right) = \exp\left(-\frac{\delta J(J+1)}{4k}\right).
$$ 
Under this condition, $\zeta^W_0$ is (conditionnally to the enriched past) the sum of independent Bernoulli variables with parameters smaller than $1/J$.
If $k$ is large enough, this happens with probability larger than $1-\varepsilon/2$. 

\subsubsection*{Control of the parameters' sum} 
Since the $F_j$ are independent and Poisson distributed with parameter $\frac{\delta j}{2k}$, the expected value of the sum $S:=\sum_{j=1}^k j^{-1} F_j $ of the parameters is $\delta/2$, and its variance is
$$ \var S = \sum_{j=1}^k \frac{1}{j^2}\frac{\delta j}{2k} = \frac{\delta}{2k} \sum_{j=1}^k \frac{1}{j}.
$$
Hence, if $k$ is large enough, 
$$ \PP \left(|S-\delta/2|<\varepsilon_1\right) > 1-\varepsilon/2. $$

Putting things together, we have proved that with probability larger than $1-\varepsilon$ on the enriched past, the conditional distribution of the number $\zeta_0^W$ of white particles at site 0 satisfies
$$ \sum_{\ell\ge0} \left| \PP(\zeta_0^W=\ell | \mbox{enriched past})-\exp(-\delta/2)\frac{(\delta/2)^\ell}{\ell !}\right| < \varepsilon. $$
This proves Lemma~\ref{lemma:general} when $\AA$ is reduced to a singleton.
\end{proof}

\begin{proof}[Proof of Lemma~\ref{lemma:general}, general case]
We consider the family of independent processes 
$\zeta^{a,b}=\left(\zeta^{B,a,b}, \zeta^{W,a,b}\right)$, $(a,b)\in \AA\times\AA$, 
which counts the number of black and white $\zeta$ particles at $x$ belonging to a pair of black and white particles respectively labelled by $a$ and $b$. 
Then $\zeta^{a,b}$ is a simple-case $\zeta$ process, for which the expected number of black particles per site is $\delta  \lambda_{\AA\times\AA}(a,b)/2$. 
From the proof in the simple case, we know that as soon as $k$ is large enough, the $\dbar$-distance between $\zeta^{a,b}$ and $\xi^{a,b}$ is smaller than $\varepsilon/|\AA|^2$, 
where $\xi^{a,b}=\left(\xi^{B,a,b}, \xi^{W,a,b}\right)$ is composed of two i.i.d. sequences of Poisson random variables of parameter $\delta  \lambda_{\AA\times\AA}(a, b)/2$.
We can recover $\zeta^B_x(a)$ and $\zeta^W_x(b)$ by
$$
\zeta^B_x(a) = \sum_{b\in\AA} \zeta_x^{B,a,b} \quad\mbox{ and }\quad\zeta^W_x(b) = \sum_{a\in\AA} \zeta_x^{W,a,b}. $$
On the other hand, $\sum_{b\in\AA} \xi^{B,a,b}$ (respectively $\sum_{a\in\AA} \xi^{W,a,b}$) has the same distribution as $\xi^B(a)$ (respectively $\xi^{W}(b)$): It is an i.i.d. sequence of Poisson random variables of parameter $P^B(a)\delta /2$ (respectively $P^W(b)\delta /2$).
Summing over $a$ and $b$, it follows that $\dbar\left(\zeta,\xi\right)<\varepsilon$.
\end{proof}

\section{Positive Poisson entropy}
\label{Sec:positiveEntropy}
We denote by $\xi^{(\infty)}$ the stationary process living in the Poisson suspension of our transformation $T$, defined by 
$$ \xi^{(\infty)}_x := \mbox{number of particles in $A$ at time }x. $$
The purpose of this section is to show that the entropy of the process $\xi^{(\infty)}$ is positive as soon as the $k_n$'s are chosen large enough. This will be proved by showing that the $\dbar$-distance between $\xi^{(\infty)}$ and an i.i.d. sequence $\xi^{(0)}$ of random Poisson variables with parameter 1 can be made as small as we want. By Lemma~\ref{lemma:entropy}, this will be enough to conclude.

Our strategy is the following: As one goes along in the construction of the return time to $A$, we define a sequence $\left(T^{(n)}\right)$ of infinite-measure-preserving transformations, approximating the final transformation $T$. Then we consider the process $\xi^{(n)}$ living in the Poisson suspension of $T^{(n)}$:
$$ \xi^{(n)}_x := \mbox{number of particles in $A$ at time $x$ for the Poisson suspension over }T^{(n)}. $$

The transformation $T^{(0)}$ is constructed by stacking infinitely many pairwise disjoints intervals of length~1, $A$ being one of them, into a doubly infinite tower, each interval being mapped onto the one just above. Therefore, as mentioned previously, $\xi^{(0)}$ is an i.i.d. sequence of Poisson variables with parameter 1. At step $n$ of the construction, we define the return time to $A$ on the subset $A_n$. This return time to $A$ will be the same for all transformations $T^{(m)}$, $m\ge n$, and for the final transformation $T$. 
By choosing the return time adequately, we will make sure that 
$$ \dbar\left(\xi^{(n-1)},\xi^{(n)}\right) < 2^{-n} \varepsilon,$$
so that for all $n$
\begin{equation}
\label{dbar-xi}
 \dbar\left(\xi^{(0)},\xi^{(n)}\right) < \varepsilon.
\end{equation}

Let us describe the first step. Recall that Tower~1 is of height 2, with basis $A_1$ (See Figure~\ref{Fig:odometer}).
We cut $A_1$ into $k_1$ equal subintervals, and we define on $A_1$ the return time to $A$ to be $M_1+j-1$ on the $j$-th subinterval. 
We insert Tower~1 into a doubly infinite tower of intervals of length 1/2 and add spacers between $A_1$ and its image by the odometer $S$: We insert $M_1+j-2$ spacers of width $1/(2k_1)$ between the $j$-th subinterval of $A_1$ and its image by $S$. 
The transformation $T^{(1)}$ is defined by mapping each point to the point right above it.

The process $\xi^{(1)}$ counts the number of particles in $A$ at time $x$ for the Poisson suspension over $T^{(1)}$:
We interpret points in $A_1$ as black particles and points in $A\setminus A_1$ as white particles. 
For the suspension over $T^{(0)}$, black and white particles are independent, whereas for the suspension over $T^{(1)}$, they are linked through the return time to $A$ of the point corresponding to the black particle.
Hence, a direct application of Lemma~\ref{lemma:general} in the simple case with no alphabet gives that 
$\dbar\left(\xi^{(0)},\xi^{(1)}\right) <  \varepsilon/2$ when $k_1$ is large enough.

\begin{center}
\begin{figure}
 \input{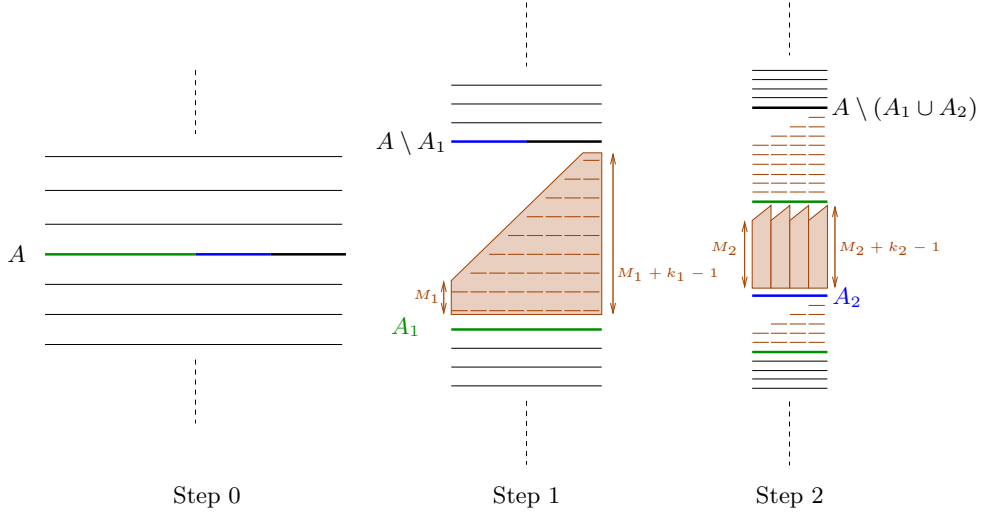}
\caption{First steps in the construction of the return time to $A$.}
 \label{Fig:induction}
\end{figure}
\end{center}

Suppose the return time to $A$ has already been defined on all $A_i$, $1\le i\le n-1$.
Consider Tower~$n$. The return time to $A$ has already been defined on all rungs but the roof and $A_n$ (which is the rung of level $2^{n-1}$). 
For each point $y$ in $A_n$, let $h^B(y)$ be the sequence of the return times to $A$ when we climb the first half of Tower~$n$ before reaching $y$, and let $h^W(y)$ be the sequence of the return times to $A$ when we climb the second half of Tower~$n$ starting from $Sy$.
We denote by $r(y)$ the return time to $A$, which is to be defined at this step. 
We want $r(y)$ to be independent of $h^B(y)$ and $h^W(y)$. 
To this end, we consider the finite partition of $A_n$ generated by $h^B$ and $h^W$. Each atom of this partition is cut into $k_n$ equal pieces, and we define $r$ to be $M_n+j-1$ on the $j$-th piece of each atom. 
Here is how we define the transformation $T^{(n)}$: we insert Tower $n$ into a doubly infinite tower of intervals of length $1/2^n$ and insert as many spacers as we need between the rungs of Tower $n$ to achieve the already defined return time to $A$. The transformation $T^{(n)}$ maps each point to the point right above it.

Let us turn to the estimation of $\dbar\left(\xi^{(n-1)},\xi^{(n)}\right)$ for $n\ge 2$.
We want to apply Lemma~\ref{lemma:general}: Black particles are points in $A_n$ and white particles are points in $SA_n$. 
Let $R_{n-1}$ be the maximum value of the already defined return time to $A$ on $A_1\cup\cdots\cup A_{n-1}$. 
We consider the finite alphabet $\AA:=\left\{1,\cdots, R_{n-1}\right\}^{2^{n-1}-1}$.
To each point $y$ in $A_n$, we associate the label $h^B(y)\in\AA$ of the return times to $A$ when we climb the first half of Tower~$n$ before reaching $y$. 
To each point $Sy$ in $SA_n$, we attach the label $h^W(y)\in\AA$, which is the sequence of the return times to $A$ when we climb the second half of Tower~$n$ starting from $Sy$.
Let the process $\xi$ (respectively $\zeta$) count the numbers of black and white particles together with their label in the suspension over $T^{(n-1)}$ (respectively $T^{(n)}$). These processes are exactly of the form studied in Lemma~\ref{lemma:general}.

Now, observe that we can recover $\xi^{(n-1)}$ and $\xi^{(n)}$ from $\xi$ and $\zeta$:
\begin{multline*}
 \xi^{(n-1)}_0 =  \sum_{h\in\AA} \xi_0^{B}(h) + \xi_{h_1}^{B}(h) + \dots +\xi_{h_1+\dots+h_{2^{n-1}-1}}^{B}(h) \\
+\sum_{h\in\AA} \xi_0^{W}(h) + \xi_{-h_1}^{W}(h) + \dots +\xi_{-(h_1+\dots+h_{2^{n-1}-1})}^{W}(h),
\end{multline*}
and
\begin{multline*}
 \xi^{(n)}_0 = \sum_{h\in\AA} \zeta_0^{B}(h) + \zeta_{h_1}^{B}(h) + \dots +\zeta_{h_1+\dots+h_{2^{n-1}-1}}^{B}(h) \\
 +\sum_{h\in\AA} \zeta_0^{W}(h) + \zeta_{-h_1}^{W}(h) + \dots +\zeta_{-(h_1+\dots+h_{2^{n-1}-1})}^{W}(h).
\end{multline*}
It follows that, if $\xi$ and $\zeta$ coincide on $\{-R_{n-1}(2^{n-1}-1),\ldots,R_{n-1}(2^{n-1}-1)\}$, then $\xi^{(n-1)}_0=\xi^{(n)}_0$. By Lemma~\ref{lemma:general}, $\dbar\left(\xi,\zeta\right)$ can be made arbitrarily small by choosing $k_n$ large enough. Hence we can assure that 
$ \dbar\left(\xi^{(n-1)},\xi^{(n)}\right) < 2^{-n} \varepsilon$. 

Finally, note that since $M_n$ is increasing, the return time to $A$ on the roof of Tower~$n$ (the union of $A_i$, $i>n$) will be larger than $M_{n+1}$. Hence, for any $L>0$, if $n$ is large enough so that $M_{n+1}>L$, the distribution of $\xi^{(n)}\big|_0^{L-1}$ coincides with the distribution of $\xi^{(\infty)}\big|_0^{L-1}$.
Therefore, 
$$ \dbar_L\left(\xi^{(\infty)},\xi^{(0)}\right) = \dbar_L\left(\xi^{(n)},\xi^{(0)}\right) \le \dbar\left(\xi^{(n)},\xi^{(0)}\right) < \varepsilon, $$
which implies 
$$ \dbar\left(\xi^{(\infty)},\xi^{(0)}\right)\le \varepsilon. $$
%

\section{Comments and open questions}

In view of previously known results on the subject, some comments on the infinite-measure-preserving transformation $T$ constructed in Section~\ref{Sec:construction} may be made. 

First, although its construction is derived from the standard cutting-and-stacking procedure used to build the most elementary rank-one system (the Von Neumann-Kakutani odometer), the transformation $T$ is not even of finite rank. Indeed, Proposition~10.1 in~\cite{jmrr08} shows that for finite-rank systems, both the Poisson and Krengel entropy vanish. 

Second, it was also proved in~\cite{jmrr08} that Krengel and Poisson entropies coincide for \emph{quasi-finite} transformations, namely transformations for which there exists a sweep-out set $A$ of measure~1 such that the return-time partition of $A$ has finite entropy. There exist only few examples of transformation for which the non quasi-finiteness has been established: an unpublished example constructed by Ornstein has been mentioned by Krengel in~\cite{krengel1969}, and the only published example which we are aware of is a rank-one system, given by Aaronson and Park in~\cite{aaronson_park2009}. Our construction thus provides a new example of a non quasi-finite transformation. 

\bigskip

After having proved that the different notions of entropy for infinite-measure-preserving transformations do not always coincide, a natural question is to ask whether they are always ordered in the same way: Is it true that Poisson entropy always dominates Krengel entropy? Can we at least decide whether zero Poisson entropy implies zero Krengel entropy? And what about similar questions regarding the comparison between Parry entropy and Poisson entropy? It may be worth recalling here that the equality of Parry entropy and Krengel entropy in the quasi-finite case was proved by Parry in 1969~\cite{parry_generators69}, but that the question whether they always coincide is, as far as we know, still open.

\subsection*{Acknowledgements}

The construction of the transformation $T$ in Section~\ref{Sec:construction} has been inspired by a private communication of Benjamin Weiss concerning the Shannon-McMillan-Breiman theorem. 

We are also much indebted to Emmanuel Roy for stimulating conversations on the subject.

\bibliography{poisson.bib}
\end{document}